\documentclass[conference,a4paper]{IEEEtran}

\IEEEoverridecommandlockouts                        

\IEEEpubid{\begin{minipage}{\textwidth}\ \normalsize \\ \\ \\ \\ \\ \\ \\ \centering Preprint submitted to 15th IEEE International Conference on Environment and Electrical Engineering (EEEIC 2015).\end{minipage}} 

\usepackage{graphicx}
\usepackage{amsmath} 
\usepackage{amssymb}
\usepackage{eurosym} 
\usepackage{booktabs}
\usepackage[switch]{lineno}
\usepackage{xcolor}
\usepackage{acronym}
\usepackage[noadjust]{cite} 

\newcommand{\bdm}{\begin{displaymath}}
\newcommand{\edm}{\end{displaymath}}
\newcommand{\bi}{\begin{itemize}}
\newcommand{\ei}{\end{itemize}}

\newcommand{\bc}{\begin{center}}
\newcommand{\ec}{\end{center}}
\newcommand{\be}{\begin{equation}}
\newcommand{\ee}{\end{equation}}
\newcommand{\bma}{\begin{math}}
\newcommand{\ema}{\end{math}}
\newcommand{\bea}{\begin{eqnarray}}
\newcommand{\eea}{\end{eqnarray}}
\newcommand{\ba}{\begin{align}}
\newcommand{\ea}{\end{align}}
\newcommand{\bal}{\begin{aligned}}
\newcommand{\eal}{\end{aligned}}
\newcommand{\barr}{\begin{array}}
\newcommand{\earr}{\end{array}}

\begin{document}
\title{Electric Energy Storage Systems integration \\in Distribution Grids}

\author{
\IEEEauthorblockN{Alessandro Di Giorgio, Francesco Liberati, Andrea Lanna}
\IEEEauthorblockA{Department of Computer, Control and  Management Engineering\\
"Sapienza" University of Rome\\
Via Ariosto 25, 00185, Rome, Italy\\
\{digiorgio, liberati, lanna\}@diag.uniroma1.it}

\thanks{This work is partially supported by the European Union FP7 ICT-GC MOBINCITY project, grant agreement no. 314328, and the SAPIENZA - ATENEO 2013 "Planning and control of flexible electricity demand and generation from renewable energy sources in Smart Grids" project, no. C26A13LYTB.}
}

\setlength{\IEEEilabelindent}{0\parindent}
\setlength{\IEEEelabelindent}{0\parindent}


\acrodef{DER}{Distributed Energy Resources}
\acrodef{DFIG}{Doubly Fed Induction Generator}
\acrodef{DSO}{Distribution System Operator}
\acrodef{ESS}{Energy Storage System}
\acrodef{ISO}{Independent System Operator}
\acrodef{MPC}{Model Predictive Control }
\acrodef{OPF}{Optimal Power Flow}
\acrodef{PLC}{Power Line Communications}
\acrodef{QP}{Quadratic Programming}
\acrodef{RES}{Renewable Energy Sources}
\acrodef{RT-OPF}{Real Time Optimal Power Flow}
\acrodef{TSO}{Transmission System Operator}
\acrodef{VPP}{Virtual Power Plant}

\maketitle

\begin{abstract}
This paper presents a real time control strategy for dynamically balancing electric demand and supply at local level, in a scenario characterized by a HV/MV substation with the presence of renewable energy sources in the form of photovoltaic generators and an electric energy storage system. The substation is connected to the grid and is powered by an equivalent traditional power plant playing the role of the bulk power system.
A Model Predictive Control based approach is proposed, by which the active power setpoints for the traditional power plant and the storage are continually updated over the time, depending on generation costs, storage's state of charge, foreseen demand and production from renewables. The proposed approach is validated on a simulation basis, showing its effectiveness in managing fluctuations of network demand and photovoltaic generation in test and real conditions.
\end{abstract}


\section*{Nomenclature} \label{Nomenclature}

\resizebox{\columnwidth}{!}{ \begin{tabular}{ c l } 
  $P^l$   			& Bus active power demand \\
  $P^g$   			& Active power generation from traditional power plant\\
  $P^s$   			& Storage active power flow\\
  $P^{res}$		& Aggregated active power of photovoltaic generators\\
  $C$			& Cost function associated to traditional power plant\\
  $x$            		& Storage's state of charge\\
  $\Theta$       		& Sampling time \\
  $T$            		& Set of time slots in the control horizon
\end{tabular}}


\section{Introduction}\label{sez:chapter1}

Unlike other types of energy, for many years it has not been possible to conveniently store electricity in significant amounts. This is the reason why a fundamental requirement for traditional power systems has been the one of operating power production so as to instantaneously meet the continuously changing load demand for active and reactive power. In this regard, adequate spinning reserve has to be maintained and controlled at all the times to assure a reliable operation of the network \cite{MACHOWSKI_BOOK11}.

As challenging targets in power production from \ac{RES} have been established at international levels \cite{SEC2007}, the share of fluctuating \ac{DER} is rapidly increasing in the advanced countries, so that a significant need for unconventional balancing power and storage capabilities aimed at progressively relaxing the above-mentioned requirement has emerged in the last years. Consequently, technological advancements in the field of electric \acp{ESS} have been obtained, putting the basis for new operational strategies of electricity systems at high, medium and low voltage levels. \IEEEpubidadjcol

Some works are becoming to appear in the relevant literature regarding this new research topic, dedicated both to theoretical, applicative and economic aspects. \ac{ESS} integration at residential level is investigated in \cite{CASTILLO_SE11,ADG_MED12b,ArboleyaTSG,LIU_PWC14,MARTIRANO_EEEIC14}, in relation to the problem of supporting consumers in the participation to demand side management programs, also taking advantage of ICT functionalities offered by Future Internet technology \cite{CASTRUCCI_LNCS11,SURACI_LNCS13,AL_IEVC14}. Other works consider higher size of ESS for applications in microgrids and distribution electricity grids \cite{FALVO_AEIT14,GILL_TPS14,LEVRON_TPS13,CERVONE_MEDPOWER12}. Also, a significant attention of academics is going towards the implications of storage integration at transmission level. Dynamic formulations of the optimal power flow problem at transmission level have been introduced and analyzed, both on a simulation basis (see e.g. \cite{WANG_PES13}) and from a theoretical perspective (see e.g. \cite{CHANDY_CDC10}, considering a simplified network constituted by a single generator and a single load, and \cite{GAYME_TPS13}, where a proof of concept is provided taking the IEEE 14 bus system as case study). Furthermore, in \cite{NAKAYAMA_SGC13} the integration of ESS and renewable energy sources is jointly investigated and simulated.
Academic and industrial worlds are currently working on this topic, with the aim of showing the technical and economic feasibility of storage integration through national and international research projects (e.g. E-Cube \cite{PR_ECUBE} and GRID4EU \cite{PR_GRID4EU}) and pioneering demonstrator in operation (e.g \cite{NEC}). 

In this paper a \ac{MPC} strategy is presented as a potential tool for facilitating the integration of a medium/large-scale electric ESS working at HV/MV substation level, with the objective of mitigating the effect of fluctuations in the aggregated distribution grid demand and in the generation from \ac{DER}. 
The \ac{MPC} methodology has been proved to be effective in a number of industrial applications \cite{MAYNE_AUT00,ELLIS_JPC14,ADG_MED13,LZ_IEVC14,SBORDONE_EPSR15}. Its working logic appears promising also in the present context: based on a feedback from the field and a continuous re-optimization which relies on a prediction model providing the state of charge over a control period ahead in time, it realizes ESS integration through a closed loop optimal control able to manage fluctuations in the boundary conditions and mitigate the effect of model inaccuracies; as a matter of fact, the proposed method overcomes conventional ESS activation strategies relying on thresholds-crossing and schedules designed in advance.
Beyond the mitigation of the net peak power flowing in the substation, which is modeled through an equivalent traditional power plant, the optimality criteria used at each iteration of the algorithm also consider the need of minimizing the distance of the storage's state of charge from a reference value (in fact, it is important to prevent the storage from being fully charged or discharged during normal grid operation, so that it can be available to absorb potential exceeding energy from \ac{RES} and support the grid in case a contingency occurs). Simulation results are focused on showing the relevance of the proposed dynamic approach, making use of both artificially generated test signals and real aggregated demand and photovoltaic generation profiles.

The remainder of the paper is organized as follows. Section \ref{sez:scenario} presents the reference scenario in which the ESS is aimed to be integrated. Section \ref{sez:modeling} is dedicated to the formalization of the control problem. In Section \ref{sez:results} the simulation results are presented and finally, in Section \ref{sez:conclusions}, the concluding remarks are drawn.

\vspace{- 0.1 cm}

\section{Reference scenario} \label{sez:scenario}

\begin{figure}[thbp]
\centering
\includegraphics[width=\columnwidth]{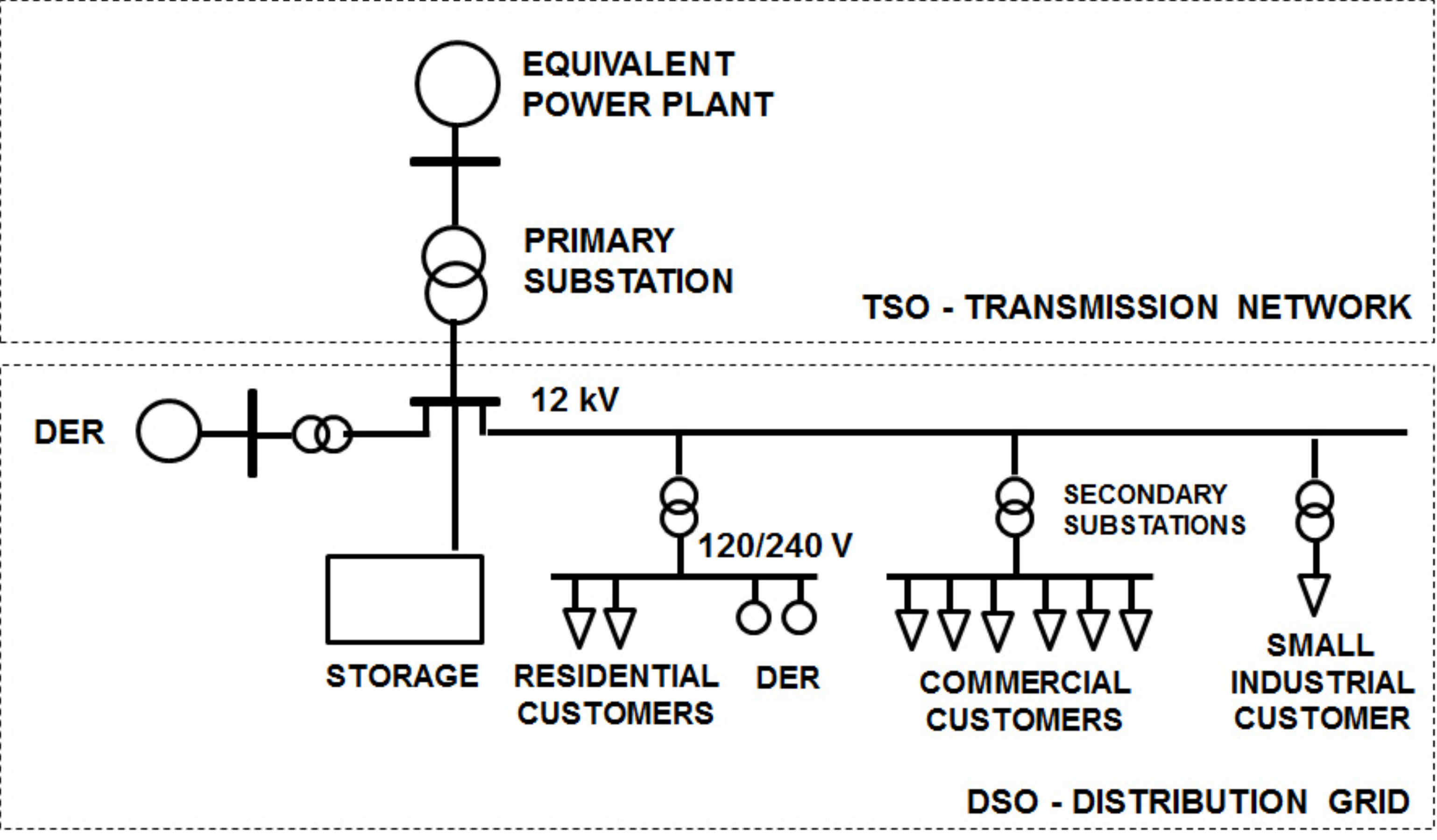}
\caption{Reference scenario.}
\label{fig:scenario}
\end{figure}

The reference scenario is depicted in Fig.~\ref{fig:scenario}. The core is represented by a HV/MV substation where \ac{DER} and different kinds of consumers are connected through the distribution lines managed by the \ac{DSO}. The net power flowing through the substation is supposed to be provided by an equivalent traditional power plant playing the role of the bulk power system. Without loss of technical generality, an electric ESS is connected to the MV busbar of the transformer in the substation, and controlled by the \ac{DSO} in order to guarantee that the net power flowing from the transmission to the distribution network remains close to a temporal power profile agreed in advance and is sufficienty smooth. This objective is fundamental to guarantee a natural extension of the traditional optimal power flow task to the case of \ac{DER} integration. As a matter of fact, in the traditional power flow problem the bus power demand is considered as a known input, due to the fact that purely passive transmission network nodes can be associated with predictable consumption patterns. Such an assumption is challenged (and hence the objective of this work becomes extremely relevant) when a significant share of \ac{DER} is integrated at distribution level, due to the fluctuating nature of the resource.

The proposed scenario is sufficiently general to cover also other similar situations, for example the one of a stand-alone local power system operated in islanding mode, in which a real power station is used to feed a large load (for example in a industrial area) with the support of \ac{RES} and an electric \ac{ESS}. Also, from the implementation point of view, the control strategy can be hosted by the SCADA control centre operated by the DSO through the communication network provided by a telco operator or directly hosted into the \ac{ESS} retrieving relevant field data making use of DSO electric infrastructure via \ac{PLC} \cite{SC_ITSM13}.

\bigskip

\section{Modeling} \label{sez:modeling}

The proposed control strategy makes use of the \ac{MPC} concept \cite{MAYNE_AUT00}. A receding control horizon is established. 
At each time period the storage's state of charge is measured and taken as initial condition for a proper dynamic prediction model. Forecasts of the bus demand and aggregated \ac{RES} generation stochastic processes are supposed to be known and updated over the control horizon.
At each MPC iteration, an open loop optimal control problem is solved, which is defined by interconnecting the balancing power problems related to all the time periods in the control horizon using the prediction model. The calculated output is represented by the power schedule for the traditional power plant and the storage device.
Reoptimization is used in order to mitigate the effects of inaccurate system modeling and take into account the updates in bus demand over the time. The control problem is formalized in the following (symbols are defined in the nomenclature section). 

The following objective function to be minimized is considered, with the aim of establishing a trade-off between the need of minimizing generation costs and that of avoiding large excursions in the storage state of charge
\be 
F = \sum_{t\in T} \left\{ \alpha(t) C(P^g(t),t) +  
    \beta(t)\left[x(t)-x^{\text{ref}}\right]^2\right\}
\label{eq:targetfunction}
\ee
The first term represents the total generation cost from thetraditional power plant over the control horizon, where each contribution can be differently weighted using the parameter $\alpha(t)$. In principle, the cost function $C(P^g(t),t)$ is a nonlinear, time variant function of the generated power. As customary in power systems studies~\cite{JIANG_SPR09,OSMAN_AMC04}, the quadratic function $C(P^g(t),t) = \gamma(t)P^{g2}(t)$ is here considered.
In the second term of the target function, $x(t)-x^{\text{ref}}$ is the excursion of the storage's state of charge with respect to its reference; this term is summed up over all the  the control horizon using different weights $\beta(t)$, in order to take bounded the evolution over all the considered time period. 

A prediction model is used to compute the storage's state of charge evolution over the control horizon. The considered dynamics is
\be
x(t) = x(t-1) - \Theta P^s(t) \qquad \forall t\in T 
\label{eq:storagedynamics}
\ee

The control $P^s(t)$ can be seen as the setpoint established by the proposed controller to drive an inner control loop locally working at storage level, equivalently to the scheduled power generation setpoint $P^g(t)$ for the prime mover governor of the traditional power plant. However we remark that differently from synchronous machines dynamics, storage dynamics is included in the proposed control problem formulation, since it is not possible to apply indefinitely the same setpoint to the storage without incurring in a saturation.

Regarding the overall behavior of the bus, the balance between demand and supply over the control horizon is included in the problem formulation through the following algebraic equation 

\be
P^l(t) - P^{res}(t) = P^g(t) + P^s(t) \qquad \forall t\in T
\label{eq:balance}
\ee

establishing that the difference between the foreseen bus power demand and generation from renewables is covered by the traditional power plant and storage.
Finally, as customary in power flow studies, proper box constraints are introduced to keep physical variables within their acceptable operating range. In this case, the restrictions to be considered are on the active power exchange and the storage's state of charge, which are limited as follows
\be
\bal
\check{P}^g &\leq{P}^g(t) \leq \hat{P}^g \qquad  &\forall t\in T\\
\check{P}^s &\leq{P}^s(t) \leq \hat{P}^s \qquad  &\forall t\in T\\
\check{x}   &\leq x(t)    \leq \hat{x}   \qquad   &\forall t\in T\\
\eal
\label{eq:box}
\ee
Finally the real time balancing problem can be stated as follows. \\
\textbf{Problem}. 
\emph{For a given power bus with known demand $P^l(t)$ and \ac{RES} generation patterns $P^{res}(t)$, generation cost function $C(P^g(t),t)$ and reference storage's state of charge $x^{ref}$, solve
\be
\min_{\bf{P^g},\bf{P^s}} F
\ee
subject to the storage dynamics (\ref{eq:storagedynamics}) with initial conditions $x^0$, and the constraints (\ref{eq:balance}) (\ref{eq:box}) over the control horizon $T$.
}

The problem previously defined can be easily put in the form of a \ac{QP} problem, which can be written as
\be
\label{QP}
\min\frac{1}{2}\textbf{u}^T\textbf{H}\textbf{u}+\textbf{f}^T\textbf{u}
\ee
\be
\label{QP1}
\textbf{A}_{ineq}\textbf{u}\leq \textbf{B}_{ineq}
\ee
\be
\label{QP2}
\textbf{A}_{eq}\textbf{u}=\textbf{B}_{eq}
\ee
\be
\label{QP3}
\textbf{u}_{min}\leq \textbf{u}\leq \textbf{u}_{max}
\ee

where vector $\textbf{u}$ of optimization variables is given by the collection of variables $P^g$ and $P^s$ at each time of problem definition. In particular, (\ref{QP}) can be derived by re-writing (\ref{eq:targetfunction}) in matrix form, after having found the explicit solutions of (\ref{eq:storagedynamics}) (i.e. after writing $x(t)$ $\forall t\in T$ in (\ref{eq:targetfunction}) as function of the current state $x$ and of the control variable $P^s$ up to time $t$).
In the same way, (\ref{eq:balance}) can be put in the form of equality constraints (\ref{QP2}), while the sets of constraints (\ref{eq:box}) can be put in the form of inequality constraints (\ref{QP1}) and box constraints~(\ref{QP3}).

\section{Results} \label{sez:results}

\subsection{Overview} 
\label{sez:overview}

The specific case study under simulation considers an HV/MV substation with two transformers, each one rated at 60 MVA. The equivalent traditional generator is characterized by unlimited nominal active power output. The ESS has $12$~MWh of battery capacity and $6$~MW of nominal power input/output. Initial conditions and control system references have been chosen as $x^0 = 0$~MWh and $x^{ref}=6$~MWh respectively. The considered sampling period is $\Theta = 5$~m (minutes), while the control horizon has been set to $2$~hours, unless differently indicated. Parameter $\gamma$ characterizing the cost function of the power plant has been set to $\gamma = 1$, while the weights appearing in the target function (\ref{eq:targetfunction}) have been chosen as independent from time as $\alpha = 1$ and $\beta = 5$.

Simulations have been performed using a Macbook 5.2, Intel Core 2 Duo, 2.13 GHz, 5GB RAM 800 Mhz DDR2 computer, running Apple OSX Yosemite (v. 10.10.1). The control framework has been built in Matlab R2012a 64 bit, and the QP problem at the basis of the MPC has been solved at each iteration by using the commercial solver IBM ILOG CPLEX v12.6.  The optimization problem to be solved at each iteration consists of 60 variables. The average detected time to complete an iteration has been $0.064$~s, including the pre-processing and post-elaboration phases.

Three test cases have been simulated for providing a proof of concept in sample conditions and using real power profiles, dedicated to the assessment of control system performance at steady state and in response to fluctuations in the bus demand and power production from \ac{DER}, also considering the duration of control horizon as testing parameter. The third test is performed using real values of load power and renewable energy generation. 

\subsection{Test case 1: fluctuation of demand} 
\label{sez:results1}

\begin{figure}[!thbp]
\centering
\includegraphics[width=\columnwidth]{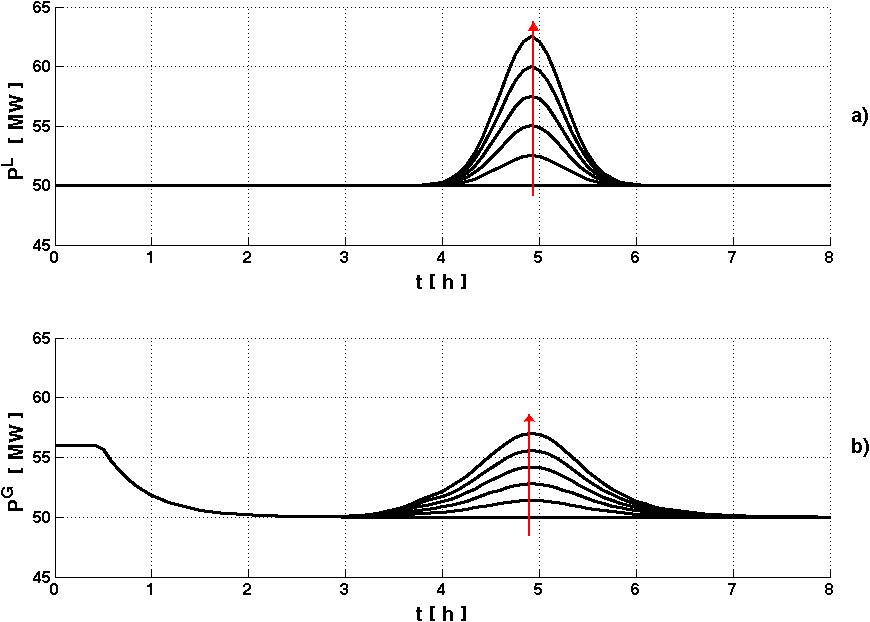}
\vspace{-12pt}
\caption{Simulation 1: (a) bus active demand, (b) active power generated by the traditional power plant.}
\label{fig:results11}
\end{figure}

\begin{figure}[!hbpt]
\centering
\includegraphics[width=\columnwidth]{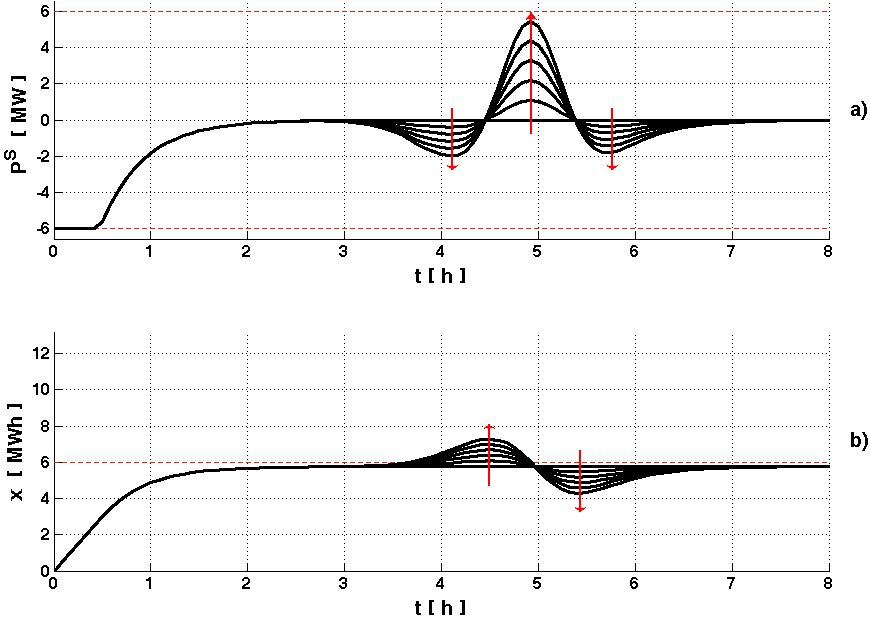}
\vspace{-12pt}
\caption{Simulation 1: storage (a) active power flow and (b) state of charge.}
\label{fig:results12}
\end{figure}

In the first simulation the foreseen aggregated demand is characterized by a base load of 50 MW and a peak is introduced by adding a Gaussian function centered at $5$~pm, having an increasing amplitude from $10\%$ to $50\%$ of the base load and $\sigma=4$, by which the deviation of the resulting network demand from the base load out of the period $4-6$~pm is less than $1 \%$ (Fig.\ref{fig:results11}a). \ac{DER} production is not included in this first simulation.
Starting from the fully discharged storage condition (see Fig.\ref{fig:results12}b, hour 0), the system then evolves towards a steady-state in which the storage state of charge is close to the reference value of $6$ MWh (the state of charge is $5.767$ MWh at hour 3, see Fig.\ref{fig:results12}b) and the network demand is completely covered by the traditional power plant (Fig.\ref{fig:results11}b, hour 3). Such an equilibrium depends on the trade-off established in the target function~(\ref{eq:targetfunction}) between generation costs and penalties for storage state of charge deviation with respect to the reference value.\par
As soon as the control horizon passes the beginning of the perturbation, the storage begins to accumulate energy, increasing the battery's state of charge (see Fig.\ref{fig:results12}). Then, as the maximum of the peak demand enters the control period, energy is released by the storage, reaching a maximum in correspondence of the maximum of bus demand (Fig.\ref{fig:results12}a). \par

After that, the storage state of charge recovers the steady-state value, ready for new interventions. It is worth noting the difference between the demand that should be covered by the traditional power plant in absence of the contribution provided by the controlled storage (Fig.\ref{fig:results11}a), and its actual power production (Fig.\ref{fig:results11}b), being the last significantly smoother than the demand. 

\begin{figure}[!th]
\centering
\includegraphics[width=\columnwidth]{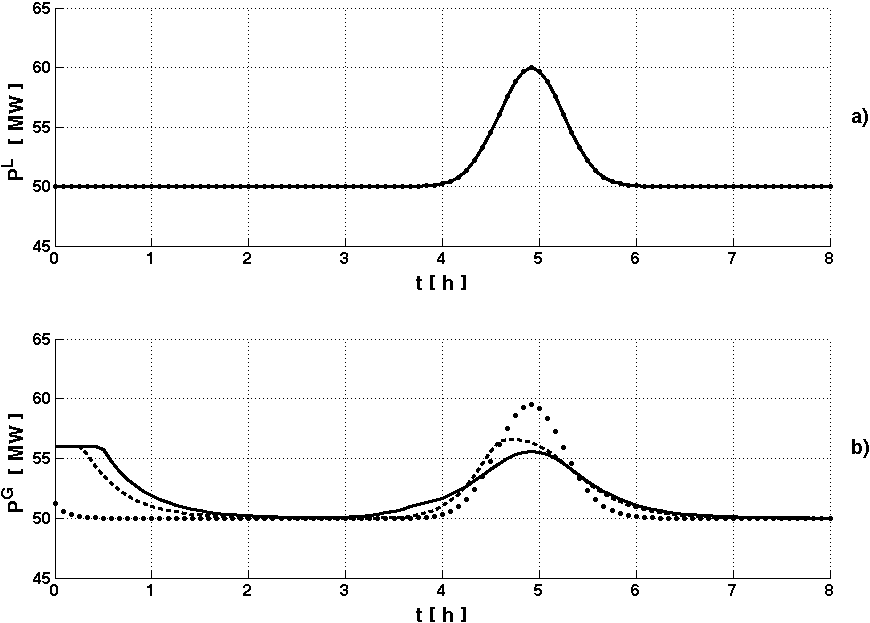}
\vspace{-12pt}
\caption{Simulation 1: (a) bus active demand, (b) active power generated by the traditional power plant when varying the control horizon.}
\vspace{-6pt}
\label{fig:results13}
\end{figure}

\begin{figure}[!th]
\centering
\includegraphics[width=\columnwidth]{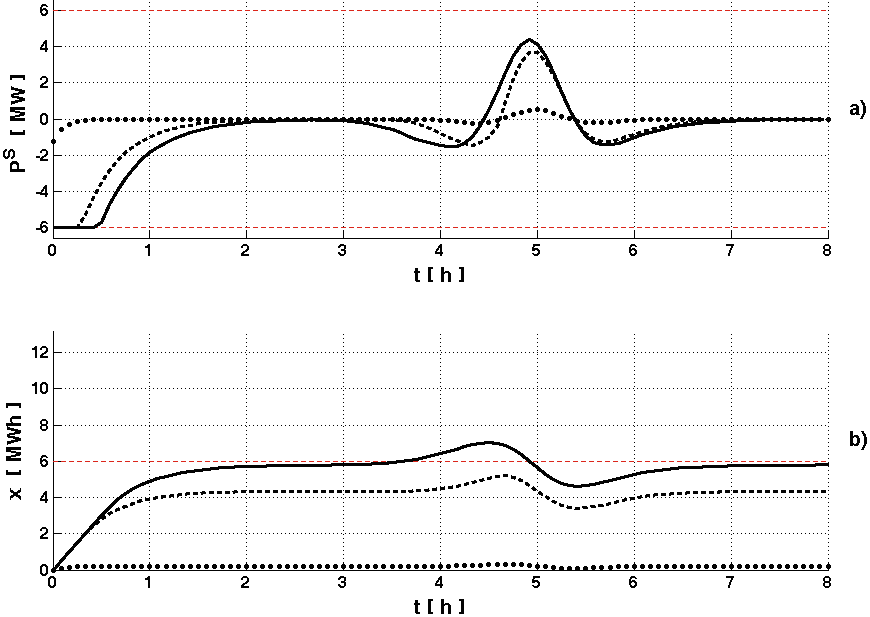}
\vspace{-12pt}
\caption{Simulation 1: storage (a) active power flow and (b) state of charge varying the control horizon.}
\vspace{-0.5 cm}
\label{fig:results14}
\end{figure}

Using the same parameters, the second part of the first simulation shows how the length of the control horizon affects the results. In Fig.\ref{fig:results13} and Fig.\ref{fig:results14}, the three curves illustrate the responses to the load profile characterized  by the peak of $20\%$ (Fig.\ref{fig:results13}a) in presence of a control horizon lasting one single time slot ($5$~m, pointed line), one hour (dotted line) and two hours (solid line). The smaller horizon has, as result, the largest distance between the storage state of charge and its reference value (Fig.\ref{fig:results14}b), which brings to a weaker peak load mitigation (Fig.\ref{fig:results13}b). Increasing the length of the control horizon, the storage accumulates energy in advance in order to intervene more effectively in presence of the maximum load (cfr. Fig.\ref{fig:results14}b). A control horizon greater than the peak size does not produce relevant improvements. \par

\subsection{Test case 2: fluctuation of RES generation} 
\label{sez:results2}

\begin{figure}[!h]
\centering
\includegraphics[width=\columnwidth]{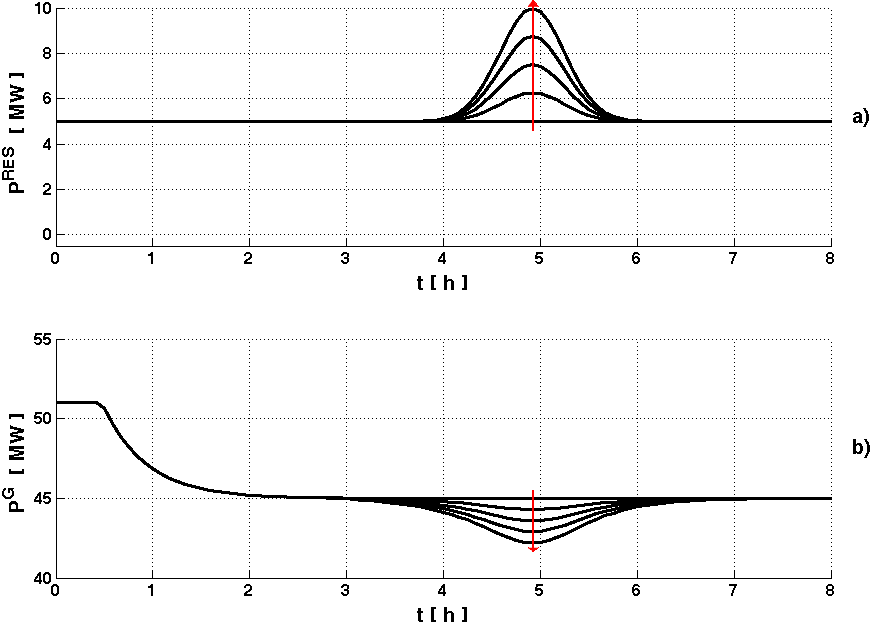}
\vspace{-12pt}
\caption{Simulation 2: (a) active power generated by DER, (b) active power generated by the traditional power plant.}
\label{fig:results21}
\end{figure}

In the second simulation, the reaction of the control system in response to a fluctuation in the power production from \ac{DER} has been tested, keeping constant the bus demand to the base value. For this purpose a base power generation of $P^{res}= 5$ MW has been considered, then an additional peak power has been introduced increasing the amplitude from $25\%$ to $100\%$ of the base value and $\sigma = 4$, by which the deviation of the resulting renewable generation from the base value out of the period $4 - 6$ pm is negligible (Fig.\ref{fig:results21}a). 

\begin{figure}[!h]
\centering
\includegraphics[width=\columnwidth]{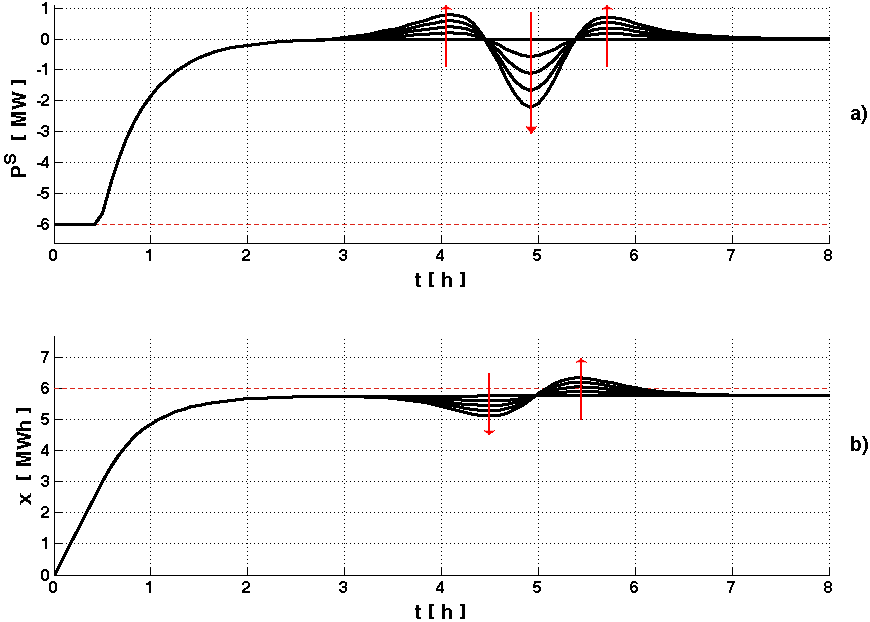}
\vspace{-12pt}
\caption{Simulation 2: storage (a) active power flow and (b) state of charge.}
\label{fig:results22}
\end{figure}

As in the previous case, once the control horizon passes the beginning of the power perturbation, the storage release energy (Fig.\ref{fig:results22}b) and then has a maximum absorption at the peak of \ac{DER} generation (Fig.\ref{fig:results22}a), then mitigating the variation of traditional plant's power production (Fig.\ref{fig:results21}b). 
\par The red arrows, similar to the previous test, illustrate the direction in which the active power produced by the \ac{RES} increases.

\subsection{Test case 3: real demand and RES generation profiles} 
\label{sez:results3}

\begin{figure}[!thbp]
\centering
\includegraphics[width=\columnwidth]{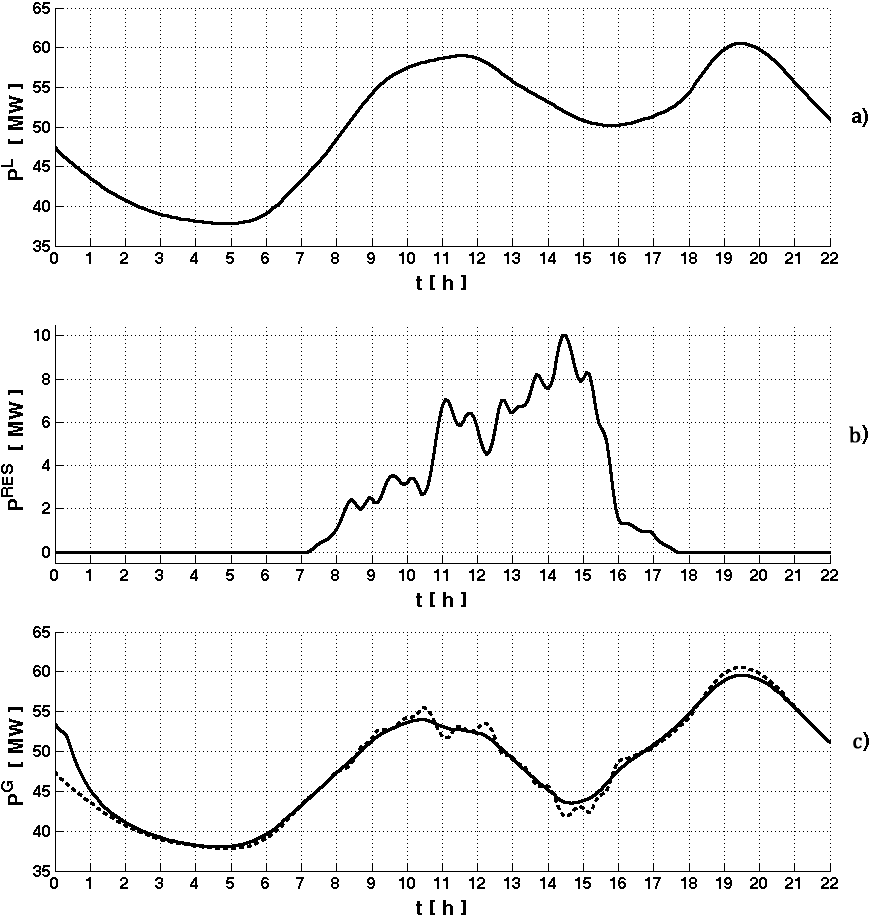}
\vspace{-12pt}
\caption{Simulation 3: (a) bus demand, (b) aggregated active power generated by photovoltaic generators, (c) active power generated by the traditional power plant.}
\label{fig:results31}
\end{figure}

\begin{figure}[!hbpt]
\centering
\includegraphics[width=\columnwidth]{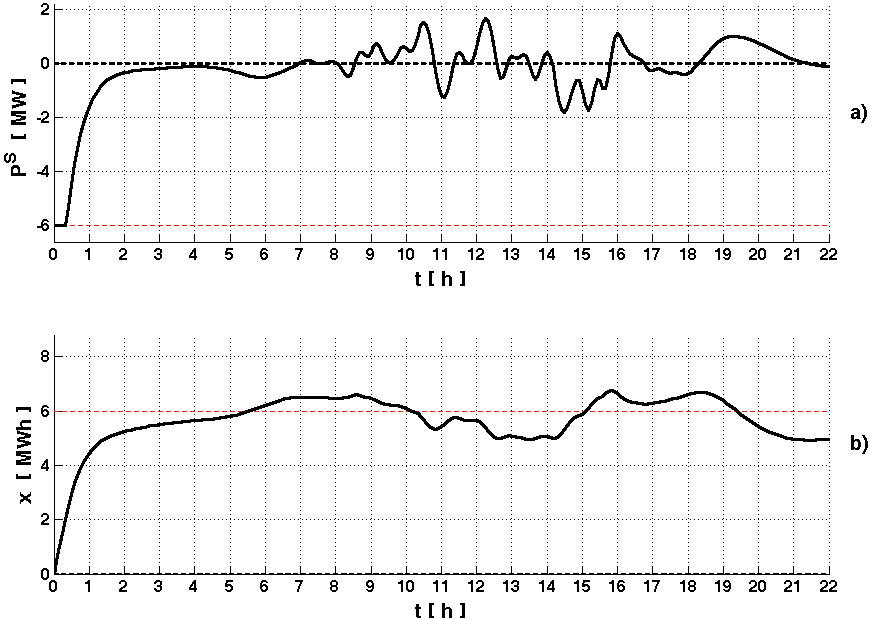}
\vspace{-12pt}
\caption{Simulation 3: storage (a) active power flow and (b) state of charge.}
\vspace{-10pt}
\label{fig:results32}
\end{figure}

In the third simulation, real bus demand (Fig.\ref{fig:results31}a, \cite{Terna}) and photovoltaic generation (Fig.\ref{fig:results31}b, \cite{PVoutput}) patterns have been considered, 
taking 2015-02-21 as tested day. Fig.\ref{fig:results31} and Fig.\ref{fig:results32} show the results of this test-case in presence of the proposed controller (solid line) and without it (dotted line).
Fig.\ref{fig:results31}c highlights the ability of the control system to mitigate fluctuations in the power production from the traditional power plant, at the cost of an initial increase of the traditional power generation in order to charge the storage (cfr. Fig.\ref{fig:results32}b). Moreover, Fig.\ref{fig:results31}c shows that the \ac{ESS} control system based on MPC, besides mitigating fluctuations of the renewable generation and, eventually, of the load, also makes smooth the most demanding power peaks (see for example between 7-9 pm).
The resulting power flow of the ESS is shown in Fig.\ref{fig:results32}a, where the variation is within $-2$ and $2$~MW. Finally, Fig.\ref{fig:results32}b shows how the the storage's state of charge remains near to the reference and tends to restore it whenever possible.
	
\section{Conclusions}\label{sez:conclusions}

A \ac{MPC} strategy has been presented as potential tool for integration of an electric \ac{ESS} in distribution electricity grids, with the objective of mitigating the effect of fluctuating aggregated demand and generation from \ac{DER} on the HV/MV substation's power flow. A proof of concept has been achieved in simplified simulation scenarios characterized by test signals and real demand/generation patterns. 
Future works will be dedicated to further testing activities aimed at assessing the robustness of the method with respect to mismatches between the forecast and the actual realization of involved stochastic processes, for example using learning technique like the one applied in \cite{ADG_CDC13} in order to better take into account the intermittent profile of renewable energy sources. Moreover, a natural evolution of the addressed problem is the introduction of power flow equations in the problem formalization for the extension of the traditional off-line optimal power flow problem to its real time dynamic version \cite{ADG_MED15}, with the aim of integrating traditionally non dispatchable generation in the operation of transmission networks. 


\section*{Acknowledgment}
The authors would like to thank Prof. Francesco Delli Priscoli, Dr. Letterio Zuccaro and Eng. Manlio Proia for the provisioning of relevant data and the helpful suggestions.




\bibliographystyle{myIEEEtran}
\bibliography{BIB_OPF,BIB_ourPaper}

\end{document}